\documentclass[final,3p]{elsarticle}
\usepackage{mathrsfs,enumerate,amssymb,amsmath,color,lineno}
\usepackage[all]{xy}
\usepackage{latexsym}
\usepackage{amsthm,a4wide,color}
\usepackage{amscd,verbatim}
\usepackage{graphicx}
\usepackage[colorlinks=true]{hyperref}

\newtheorem{theorem}{Theorem}
\newtheorem{lemma}[theorem]{Lemma}
\newtheorem{proposition}[theorem]{Proposition}

\newtheorem{example}[theorem]{Example}
\newtheorem{question}[theorem]{Question}

\newtheorem{remark}[theorem]{Remark}

\newcommand{\pf}{\noindent \textbf{Proof.}\quad}
\newcommand{\epf}{\hspace{\stretch{1}}$\blacksquare$}
\begin{document}

\begin{frontmatter}

\title{Answers to some questions about Zadeh's extension\\ on metric spaces\tnoteref{mytitlenote}}
\tnotetext[mytitlenote]{This work was supported by the National Natural Science Foundation of China (Nos. 11601449 and 11701328),
the National Nature Science Foundation of China (Key Program) (No. 51534006), Science and Technology Innovation Team of Education
Department of Sichuan for Dynamical System and its Applications (No. 18TD0013), Youth Science and Technology Innovation Team of
Southwest Petroleum University for Nonlinear Systems (No. 2017CXTD02), scientific research starting project of Southwest Petroleum
University (No. 2015QHZ029),  Shandong Provincial Natural Science Foundation, China (Grant
ZR2017QA006), and Young Scholars Program of Shandong University, Weihai (No. 2017WHWLJH09).}

\author[a1,a2]{Xinxing Wu}
\address[a1]{School of Sciences, Southwest Petroleum University, Chengdu, Sichuan 610500, China}
\address[a2]{Department of Electronic Engineering, City University of Hong Kong,
Hong Kong SAR, China}
\ead{wuxinxing5201314@163.com}

\author[a3]{Xu Zhang\corref{mycorrespondingauthor}}
\cortext[mycorrespondingauthor]{Corresponding author}
\address[a3]{Department of Mathematics, Shandong University, Weihai, Shandong 264209, China}
\ead{xu$\_$zhang$\_$sdu@mail.sdu.edu.cn}

\author[a2]{Guanrong Chen}
\ead{gchen@ee.cityu.edu.hk}

\begin{abstract}
This paper shows that there exists a contraction whose Zadeh's extension is not a contraction
under the Skorokhod metric, answering negatively Problems 5.8 and 5.12 posted in
\cite[Jard\'{o}n, S\'{a}nchez, and Sanchis, Some questions about Zadeh's extension on
metric spaces, Fuzzy Sets and Systems, 2018]{JSS2018}.
\end{abstract}


\end{frontmatter}

Zadeh's extension principle is the soul of the fuzzy set theory, and is the base for the concepts of fuzzy numbers and fuzzy arithmetic.
The Skorokhod topology is defined on the space of functions from the unit interval to the real line, where these functions are right continuous and their left limits exist. This topology is used in the study of the convergence of the probability measures, the central limit theorems and many other results in stochastic processes \cite{Bill1968,Jacod1987}.  There is a close relation between this topology and fuzzy numbers. The Skorokhod metric, induced by the Skorokhod topology, is used in the study of fuzzy numbers \cite{JSS2018,JK00}.

Throughout this paper, denote $\mathbb{N}=\left\{1, 2, 3, \ldots\right\}$ and $\mathbb{Z}^+=\left\{0, 1, 2, \ldots\right\}$. A {\it dynamical system} is a pair $(X, f)$, where $X$ is a metric space with a metric $d$ and $f: X\to X$ is a continuous map. Let ${K}(X)$ be the hyperspace on $X$, i.e., the space of non-empty
compact subsets of $X$ with the Hausdorff metric $d_{H}$ defined by
\begin{align*}
d_{H}(A, B)&=\max\left\{\max_{x\in A}\min_{y\in B}d(x, y), \max_{y\in B}\min_{x\in A}d(x, y)\right\}\\
&=\inf\left\{\varepsilon>0: A\subset B^{\varepsilon}
\text{ and } B\subset A^{\varepsilon}\right\},
\end{align*}
for $A, B\in {K}(X)$, where $A^{\varepsilon}$ is the $\varepsilon$-neighbourhood of the set $A$.

A fuzzy set $u$ in the space $X$ is a function $u: X\to I$, where $I=[0, 1]$. Given a fuzzy set $u$, its {\it $\alpha$-cuts}
(or {\it $\alpha$-level sets}) $[u]_{\alpha}$ ($\alpha\in (0, 1]$) and {\it support} $[u]_{0}$ are defined respectively by
$$
[u]_{\alpha}=u^{-1}([\alpha, 1])=\{x\in X: u(x)\geq \alpha\},
$$
and
$$
[u]_{0}=\overline{\{x\in X: u(x)>0\}}.
$$
Clearly, $[u]_{0}=\overline{\bigcup_{\alpha\in (0, 1]}[u]_{\alpha}}$.
Let $\mathfrak{F}(X)$ be the set of all upper semi-continuous fuzzy sets $u: X\to I$, satisfying that
$[u]_{0}$ is compact and $[u]_{1}\neq \O$. Define a {\it level-wise metric} $d_{\infty}$ on $\mathfrak{F}(X)$
by
\begin{equation}\label{dinf}
d_{\infty}(u, v)=\sup\{d_{H}([u]_{\alpha}, [v]_{\alpha}): \alpha\in I\}, \quad \forall u, v\in \mathfrak{F}(X).
\end{equation}

For the level-wise metric $d_{\infty}$, the following result shows that the supports are not essential for the calculation
of $d_{\infty}$.
\begin{proposition}\label{Pro-1}
For any $u, v\in \mathfrak{F}(X)$, $d_{\infty}(u, v)=\sup\{d_{H}([u]_{\alpha}, [v]_{\alpha}): \alpha\in (0, 1]\}$.
\end{proposition}

\pf
Since for any $\alpha, \beta\in (0, 1]$ with $\alpha<\beta$, $[u]_{\alpha}\supset [u]_{\beta}$ and
$[u]_{\alpha}\subset [u]_0$, it follows that $\lim_{\alpha\to 0^{+}}d_{H}([u]_{\alpha}, [u]_{0})$ exists, denoted
by $\xi$. It can be shown that $\xi=0$. In fact, if $\xi>0$, then for any $\alpha\in (0, 1]$, $d_{H}([u]_{\alpha}, [u]_{0})>\frac{\xi}{2}$.
The compactness of $[u]_0$ implies that there exist several points $x_{1}, x_{2}, \ldots, x_{n}\in [u]_{0}$ such that
$\bigcup_{i=1}^{n}B(x_{i}, \frac{\xi}{8})\supset [u]_0$, where $B(x, \varepsilon)=\{y\in X: d(y, x)<\varepsilon\}$.
Applying $[u]_{0}=\overline{\bigcup_{\alpha\in (0, 1]}[u]_{\alpha}}$ yields that for any $1\leq i\leq n$,
$B(x_{i}, \frac{\xi}{8})\cap (\bigcup_{\alpha\in (0, 1]}[u]_{\alpha})\neq \O$, i.e., there exist $\alpha_{i}\in (0, 1]$
and $z_{i}$ such that $z_{i}\in B(x_{i}, \frac{\xi}{8}) \cap [u]_{\alpha_{i}}$. Take $\alpha=\min\{\alpha_{i}: 1\leq i\leq n\}$.
Clearly, $\{z_{i}: 1\leq i\leq n\}\subset [u]_{\alpha}$. For any $x\in [u]_0$, there exists $1\leq i\leq n$ such that
$x\in B(x_{i}, \frac{\xi}{8})$, implying that $d(x, z_{i})\leq d(x, x_{i})+d(x_{i}, z_{i})<\frac{\xi}{4}$, i.e., $[u]_{0}\subset
([u]_{\alpha})^{\frac{\xi}{4}}$, where $([u]_{\alpha})^{\frac{\xi}{4}}$ is the $\frac{\xi}{4}$-neighborhood of $[u]_{\alpha}$. Thus,
$$
\xi=d_{H}([u]_0, [u]_{\alpha})\leq \frac{\xi}{4},
$$
which is a contradiction.

Let $\eta=\sup\{d_{H}([u]_{\alpha}, [v]_{\alpha}): \alpha\in (0, 1]\}$. For any $\alpha\in (0, 1]$, one has
\begin{align*}
d_{H}([u]_0, [v]_0)&\leq d_{H}([u]_{0}, [u]_{\alpha})+d_{H}([u]_{\alpha}, [v]_{\alpha})+d_{H}([v]_{\alpha}, [v]_0)\\
&\leq
d_{H}([u]_{0}, [u]_{\alpha})+ \eta + d_{H}([v]_{\alpha}, [v]_0),
\end{align*}
implying that
$$
d_{H}([u]_0, [v]_0)\leq
\eta+ \lim_{\alpha\to 0^{+}}(d_{H}([u]_{0}, [u]_{\alpha})+ d_{H}([v]_{\alpha}, [v]_0))=\eta.
$$
Therefore, $d_{\infty}(u, v)=\eta$.
\epf

\begin{remark}
{\rm In \eqref{dinf}, the value of $\alpha$ is taken from the whole interval $I=[0,1]$.}
\end{remark}

{\it Zadeh's extension} of a dynamical system $(X, f)$ is a map $\widetilde{f}: \mathfrak{F}(X)\to \mathfrak{F}(X)$
defined by
$$
\widetilde{f}(u)(x)=
\begin{cases}
0, & f^{-1}(x)= \O, \\
\sup\{u(z): z\in f^{-1}(x)\}, &  f^{-1}(x)\neq \O.
\end{cases}
$$

Let $f_1, f_{2}, \ldots, f_{n}: X\to X$ be continuous maps. Define $F: \mathfrak{F}(X)\to \mathfrak{F}(X)$ by
$$
[F(u)]_{\alpha}=[\widetilde{f}_1(u)]_{\alpha} \cup \cdots \cup [\widetilde{f}_{n}(u)]_{\alpha}, \quad \forall
u\in \mathfrak{F}(X),\ \alpha\in I.
$$

Let $\mathrm{Hom}(I)$ be the family of all strictly increasing homeomorphisms from $I$ onto itself. For any $t\in \mathrm{Hom}(I)$
and $u\in \mathfrak{F}(X)$, denote $tu=t\circ u$ for convenience. Clearly, $tu\in \mathfrak{F}(X)$.
Given a metric space $(X, d)$, the {\it Skorokhod metric}
$d_{0}$ on $\mathfrak{F}(X)$ is defined as follows \cite{JK00}:
$$
d_{0}(u, v)=\inf\left\{\varepsilon: \exists t\in \mathrm{Hom}(I) \text{ such that }
\sup\{|t(x)-x|: x\in I\}\leq \varepsilon \text{ and } d_{\infty}(u, tv)\leq \varepsilon\right\}.
$$

Recently, Jard\'{o}n et al. \cite{JSS2018} proved that both $\widetilde{f}$ and $F$ are contractive if
the previous dynamical systems are contractive under the level-wise metric $d_{\infty}$, and they proposed
the following two questions. For more recent results on Zadeh's extension and $g$-fuzzification,
refer to \cite{BK2017,FSS2018,K2016,WC2017,WDLW2017,WWL2018} and some references therein.

\begin{question}\label{Q-1}{\rm \cite[Problem~5.8]{JSS2018}
Let $(X, d)$ be a metric space and $f: X\to X$ be a contraction. Is Zadeh's extension
$\widetilde{f}: (\mathfrak{F}(X), d_0)\to (\mathfrak{F}(X), d_0)$ a contraction?}
\end{question}

\begin{question}\label{Q-2}{\rm \cite[Problem~5.12]{JSS2018}
If $(X, d)$ is a metric space and $f_1, \ldots, f_n: X\to X$ are contractions. Is
$F: (\mathfrak{F}(X), d_0)\to (\mathfrak{F}(X), d_0)$ a contraction?}
\end{question}

This paper constructs a contraction on $I$ whose Zadeh's extension is not contractive
under the Skorokhod metric $d_0$, answering negatively Questions \ref{Q-1} and \ref{Q-2} above
(see Example~\ref{Main-Exa}, Theorem \ref{Main-Thm}, and Remark~\ref{Remark-9} below).

\begin{lemma}\label{Lemma-1}{\rm\cite[Proposition 3.1]{JSS2018}}
Let $X$ be a Hausdorff space. If $f: X\to X$ is a continuous function, then for any $u\in \mathfrak{F}(X)$
and any $\alpha\in I$, one has
\begin{enumerate}[(1)]
\item $[\widetilde{f}(u)]_{\alpha}=f([u]_{\alpha})$;
\item $[tu]_{\alpha}=[u]_{t^{-1}(\alpha)}$ for $t\in \mathrm{Hom}(I)$.
\end{enumerate}
\end{lemma}

\begin{lemma}\label{Lemma-2}
Let $X$ be a metric space and $f: X\to X$ be a contraction. Then, for any $u, v\in \mathfrak{F}(X)$,
$d_{0}(\widetilde{f}(u), \widetilde{f}(v))\leq d_{0}(u, v)$.
\end{lemma}
\pf
Let $\lambda\in [0, 1)$ be a contraction factor of $f$. Applying \cite[Proposition~5.7]{JSS2018},
one has that for any $t\in \mathrm{Hom}(I)$,
\begin{align*}
d_{\infty}(\widetilde{f}(u), t\widetilde{f}(v))
&=\sup\{d_H([\widetilde{f}(u)]_{\alpha}, [t\widetilde{f}(v)]_{\alpha}): \alpha\in I\}\\
&=\sup\{d_H(f([u]_{\alpha}), f([v]_{t^{-1}(\alpha)})): \alpha\in I\}\\
&\leq \sup\{\lambda \cdot d_H([u]_{\alpha}, [v]_{t^{-1}(\alpha)}): \alpha\in I\}\\
&= \lambda\cdot d_{\infty}(u, tv).
\end{align*}
This implies that $d_{0}(\widetilde{f}(u), \widetilde{f}(v))\leq d_{0}(u, v)$.
\epf

\begin{example}\label{Main-Exa}
{\rm Construct a function $t: I\to I$ as follows:
$$
t(x)=
\begin{cases}
\frac{a}{a-\frac{1}{4}}x, & x\in [0, a-\frac{1}{4}], \\
\frac{1}{2}(x-a+\frac{1}{4})+a, & x\in [a-\frac{1}{4}, a+\frac{1}{4}],\\
x, & x\in [a+\frac{1}{4}, 1].
\end{cases}
$$
Take $a=\frac{3}{8}$, so that
$$
t(x)=
\begin{cases}
3x, & x\in [0, \frac{1}{8}], \\
\frac{1}{2}x+\frac{5}{16}, & x\in [\frac{1}{8}, \frac{5}{8}],\\
x, & x\in [\frac{5}{8}, 1].
\end{cases}
$$
\begin{figure}[h]
\begin{center}
\scalebox{0.5 }{ \includegraphics{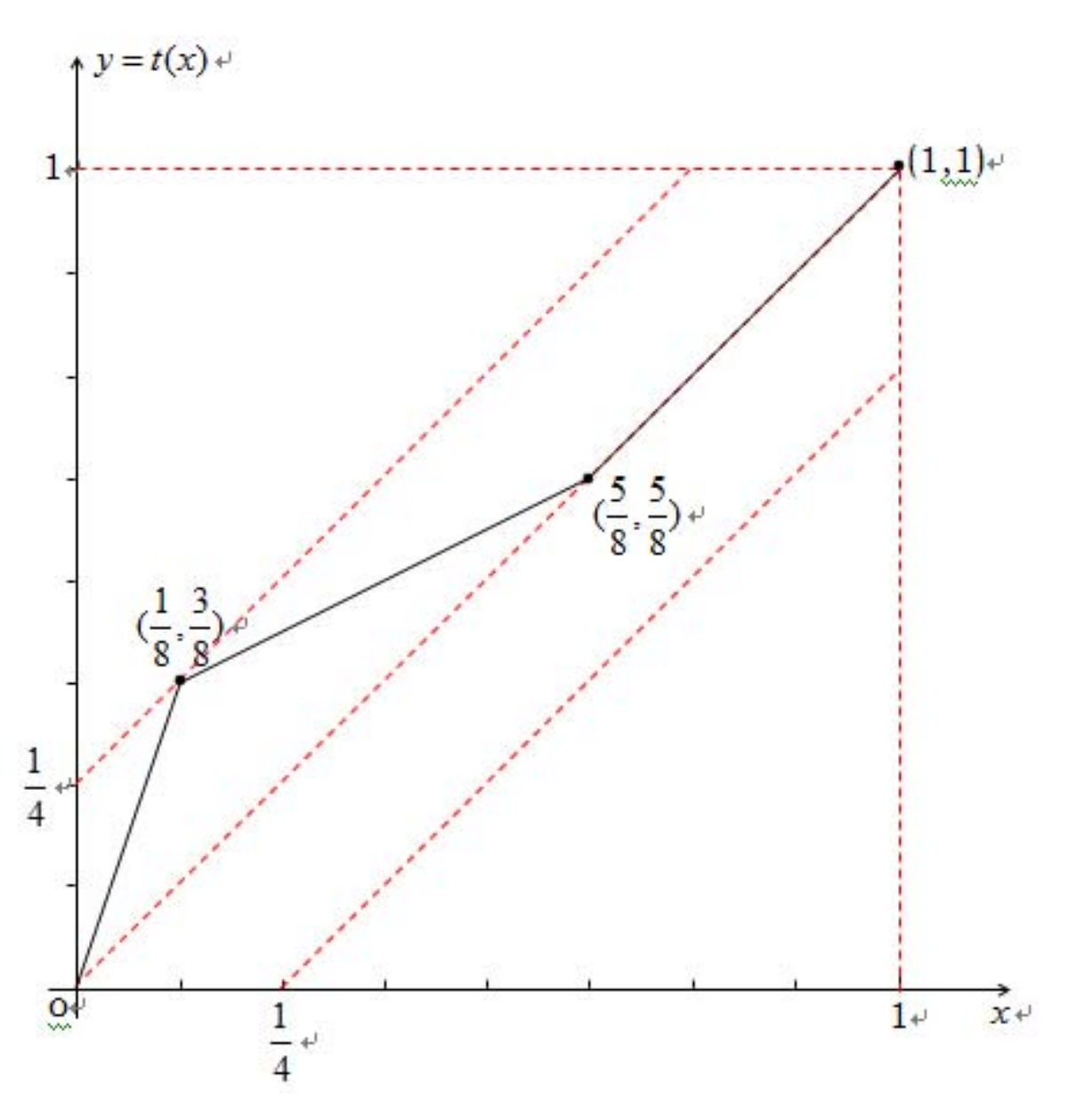}}
\renewcommand{\figure}{Fig.}
\caption{The illustration diagram of the function $t(x)$
}
\end{center}
\end{figure}
It is easy to see that $t\in \mathrm{Hom}(I)$ and
\begin{equation}\label{e-1}
\sup\left\{|t(x)-x|: x\in I\right\}=\frac{1}{4}.
\end{equation}

Take $a_0=a$, $a_1=a_{0}+\frac{1-a}{2}$, $\ldots$, $a_{n+1}=a_{n}+\frac{1-a}{2^{n+1}}$, and $b_{n}=t^{-1}(a_{n})$
for all $n\in \mathbb{N}$, and choose $a'_{0}=a$, $a'_{1}=a+\frac{1}{2}$, $a'_{2}=a'_{1}+\frac{1-(a+\frac{1}{2})}{2}$,
$\ldots$, $a'_{n+1}=a'_{n}+\frac{1-(a+\frac{1}{2})}{2^n}$ for all $n\geq 2$.
It can be verified that
\begin{enumerate}[(a)]
\item $b_{n}=a_{n}$ for all $n\geq 1$;
\item\label{i-i} $b_{1}=\frac{11}{16}>\frac{5}{8}=a+\frac{1}{4}$;
\item $\lim_{n\to +\infty}a_n=\lim_{n\to +\infty}a'_n=1$.
\end{enumerate}
Meanwhile, take $b'_{0}=a-\frac{1}{4}=\frac{1}{8}$, $b'_{1}=a+\frac{1}{2}=\frac{7}{8}$, $b'_{2}=b'_{1}+\frac{1-(a+\frac{1}{2})}{2}$,
$\ldots$, $b'_{n+1}=b'_{n}+\frac{1-(a+\frac{1}{2})}{2^n}$ for all $n\geq 2$.
Define $u: I\to I$ and $v: I\to I$ by
$$
u(x)=
\begin{cases}
0, & x\in [0, a'_0), \\
a_n, & x\in [a'_n, a'_{n+1}) \text{ and } n\in \mathbb{Z}^+,\\
1, & x=1,
\end{cases}
$$
and
$$
v(x)=
\begin{cases}
0, & x\in [0, b'_0), \\
a-\frac{1}{4}, & x\in [b'_0, b'_1), \\
a+\frac{1}{4}, & x\in [b'_1, b'_{2}),\\
b_{n-1}, & x\in [b'_n, b'_{n+1}) \text{ and } n \geq 2,\\
1, & x=1,
\end{cases}
$$
respectively.

\begin{figure}[h]
\begin{center}
\scalebox{0.5 }{ \includegraphics{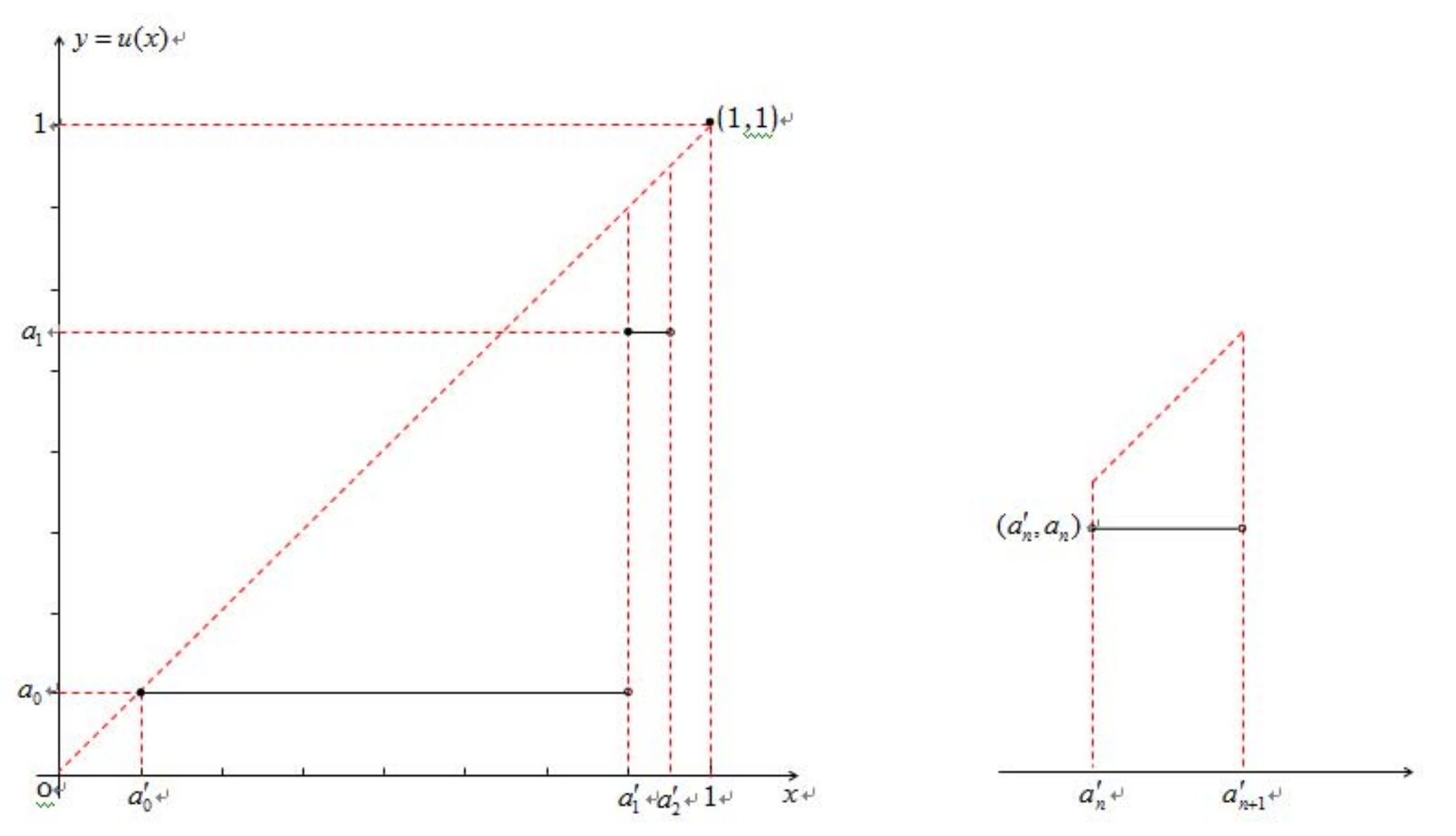}}
\renewcommand{\figure}{Fig.}
\caption{The illustration diagram of the construction of the function $u(x)$
}
\end{center}
\end{figure}

\begin{figure}[h]
\begin{center}
\scalebox{0.5 }{ \includegraphics{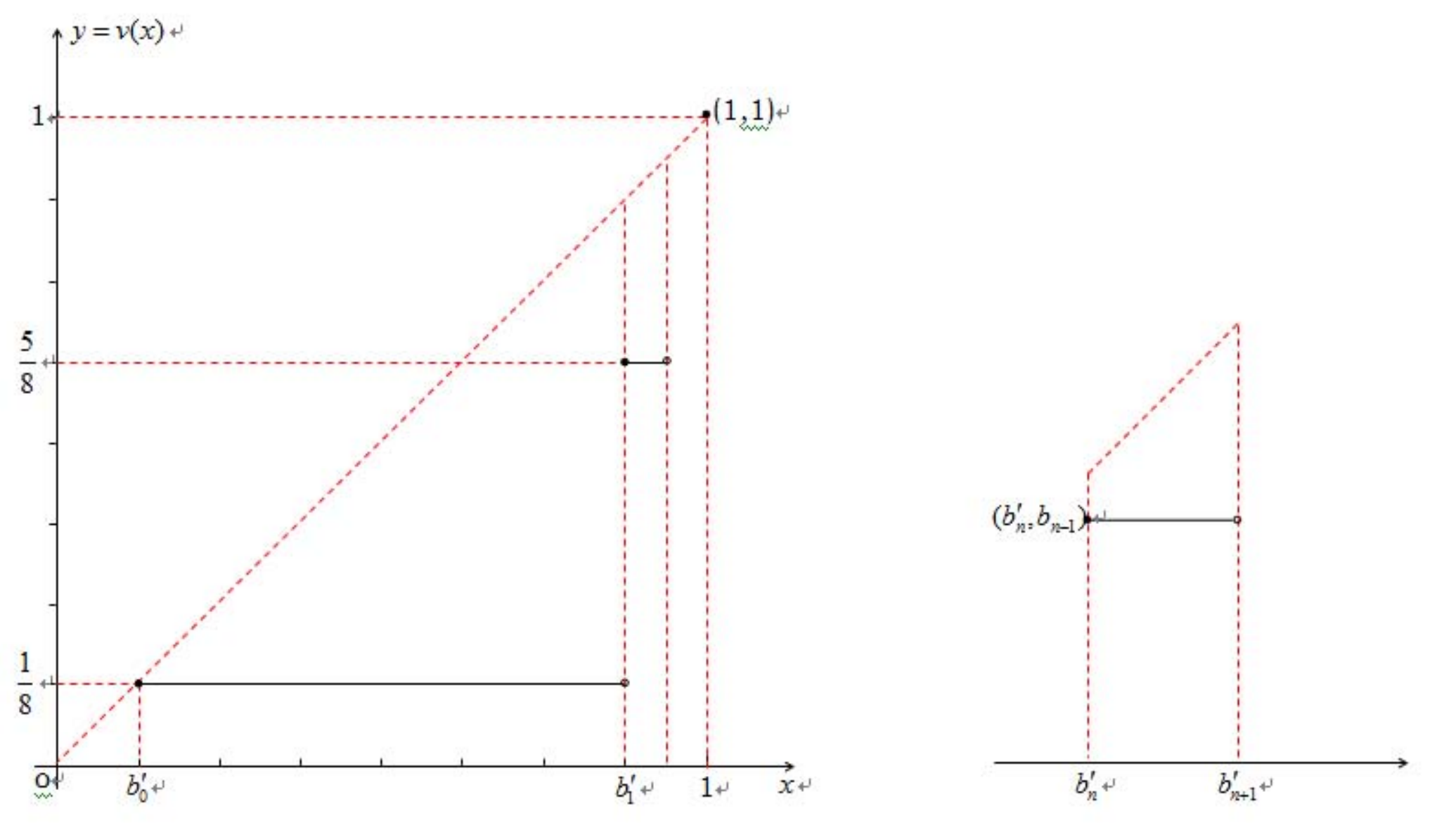}}
\renewcommand{\figure}{Fig.}
\caption{The illustration diagram of the construction of the function $v(x)$
}
\end{center}
\end{figure}

Clearly, $u, v\in \mathfrak{F}(I)$ and $u, v$ are monotonically increasing (applying (\ref{i-i})).
According to the constructions of $u$ and $v$, it can be verified that
\begin{enumerate}[(i)]
\item\label{a} $u^{-1}(1)=v^{-1}(1)=\{1\}$;
\item\label{b} $[u]_{0}=[a, 1]$;
\item\label{c} $[v]_{0}=[a-\frac{1}{4}, 1]$;
\item\label{d} $u^{-1}([x, 1])=[a, 1]$  for any $x\in (0, a]$;
\item\label{e} $u^{-1}([x, 1])=[a'_{n+1}, 1]$ for any $x\in (a_n, a_{n+1}]$, $n\in \mathbb{Z}^+$;
\item\label{f} $v^{-1}([x, 1])=[a-\frac{1}{4}, 1]$ for any $x\in (0, a-\frac{1}{4}]$;
\item\label{g}  $v^{-1}([x, 1])=[a+\frac{1}{2}, 1]$ for any $x\in (a-\frac{1}{4}, a+\frac{1}{4}]$;
\item\label{h} $v^{-1}([x, 1])=[b'_2, 1]$ for any $x\in (a+\frac{1}{4}, b_1]$;
\item\label{i} $v^{-1}([x, 1])=[b'_{n+2}, 1]$ for any $x\in (b_n, b_{n+1}]$, $n\in \mathbb{N}$.
\end{enumerate}

{\bf Claim 1.} $d_{0}(u, v)=\frac{1}{4}$.

\medskip

Fix any $0<\varepsilon<\frac{1}{4}$. For any $t'\in \mathrm{Hom}(I)$ with $\sup\{|t'(x)-x|: x\in I\}\leq \varepsilon$,
noting that $|t'(t'^{-1}(a))-t'^{-1}(a)|\leq \varepsilon$, i.e., $a-\frac{1}{4}<a-\varepsilon<t'^{-1}(a)
<a+\varepsilon<a+\frac{1}{4}$), it follows from (\ref{d}) and (\ref{g}) that
\begin{align*}
d_{\infty}(u, t'v)
&\geq d_{H}([u]_{a}, [t'v]_{a})=d_{H}([u]_{a}, [v]_{t'^{-1}(a)})\\
&= d_{H}(u^{-1}([a, 1]), v^{-1}([t'^{-1}(a), 1]))\\
&= d_{H}\left([a, 1], \left[a+\frac{1}{2}, 1\right]\right)=\frac{1}{2}.
\end{align*}
This implies that
\begin{equation}\label{e-2}
d_{0}(u, v)\geq \frac{1}{4}.
\end{equation}

(1) For any $x\in (0, a]$, from $t^{-1}(x)\in (0, a-\frac{1}{4}]$, (\ref{d}) and (\ref{f}), it follows that
$$
d_{H}([u]_{x}, [tv]_{x})=d_{H}(u^{-1}([x, 1]), v^{-1}([t^{-1}(x), 1]))=d_{H}\left([a, 1], \left[a-\frac{1}{4}, 1\right]\right)=\frac{1}{4}.
$$

(2) For any $x\in (a, a_{1}]$, from $t^{-1}(x)\in (t^{-1}(a), t^{-1}(a_1)]=(a-\frac{1}{4}, b_1]$, (\ref{e}),
(\ref{g}), and (\ref{h}), it follows that
\begin{align*}
d_{H}([u]_{x}, [tv]_{x})
&=d_{H}(u^{-1}([x, 1]), v^{-1}([t^{-1}(x), 1]))\\
&\leq \max\left\{d_{H}\left([a'_1, 1], \left[a+\frac{1}{2}, 1\right]\right), d_{H}([a'_{1}, 1], [b'_2, 1])\right\}\\
&\leq \frac{1}{4}.
\end{align*}

(3) For any $x\in (a_n, a_{n+1}]$ and any $n\in \mathbb{N}$, from $t^{-1}(x)\in (t^{-1}(a_{n}), t^{-1}(a_{n+1})]=
(b_{n}, b_{n+1}]$, (\ref{e}), and (\ref{i}), it follows that
$$
d_{H}([u]_{x}, [tv]_{x})=d_{H}(a'_{n+1}, b'_{n+2})<\frac{1}{4}.
$$

Applying Proposition~\ref{Pro-1}, one has that
$$
d_{\infty}(u, tv)=\sup\left\{d_{H}([u]_{x}, [tv]_{x}): x\in (0, 1]\right\}=\frac{1}{4}.
$$
This, together with (\ref{e-2}), implies that
$$
d_{0}(u, v)=\frac{1}{4}.
$$

For any $\lambda\in [\frac{1}{2}, 1)$, define $f_{\lambda}: I\to I$ by $f(x)=\lambda x$ for all $x\in I$.
Clearly, $f_{\lambda}$ is a contraction.

\medskip

{\bf Claim 2.} $d_{0}(\widetilde{f}_{\lambda}(u), \widetilde{f}_{\lambda}(v))=d_{0}(u, v)=\frac{1}{4}$,
and thus $\widetilde{f}_{\lambda}$ is not a contraction.

\medskip

For any $t'\in \mathrm{Hom}(I)$ with $\sup\left\{|t'(x)-x|: x\in I\right\}<\frac{1}{4}$,
from $a-\frac{1}{4}<t'^{-1}(a)<a+\frac{1}{4}$ and Lemma \ref{Lemma-1}, it follows that
\begin{align*}
d_{\infty}(\widetilde{f}_{\lambda}(u), t'\widetilde{f}_{\lambda}(v))
&\geq d_{H}([\widetilde{f}_{\lambda}(u)]_{a}, [t'\widetilde{f}_{\lambda}(v)]_{a})
=d_{H}(f_{\lambda}([u]_{a}), f_{\lambda}([v]_{t'^{-1}(a)}))\\
&=d_{H}\left(f_{\lambda}([a, 1]), f_{\lambda}\left(\left[a+\frac{1}{2}, 1\right]\right)\right)
=d_{H}\left([\lambda\cdot a, \lambda],
\left[\lambda \cdot a+\frac{\lambda}{2}, \lambda\right]\right)\\
&=\frac{\lambda}{2}\geq \frac{1}{4},
\end{align*}
implying that
$$
d_{0}(\widetilde{f}_{\lambda}(u), \widetilde{f}_{\lambda}(v))\geq \frac{1}{4}.
$$
This, together with Lemma \ref{Lemma-2}, implies that
$$
d_{0}(\widetilde{f}_{\lambda}(u), \widetilde{f}_{\lambda}(v))=\frac{1}{4}.
$$
}
\end{example}

\begin{theorem}\label{Main-Thm}
There exists a contraction $(I, f)$ such that its Zadeh's extension
$(\mathfrak{F}(I), \widetilde{f})$ is not a contraction under the Skorokhod metric $d_0$.
\end{theorem}

\begin{remark}\label{Remark-9}
{\rm \begin{enumerate}[(1)]
\item Theorem~\ref{Main-Thm} shows that the answer to Question~\ref{Q-1} is negative.
\item Choose $f_{1}, f_{2}: I\to I$ as $f_{1}=\frac{1}{2}x$ and $f_{2}(x)=\frac{3}{4}x$ for all
$x\in I$. For any $t'\in \mathrm{Hom}(I)$ with $\sup\{|t'(x)-x|: x\in I\}<\frac{1}{4}$, it can be
verified that
\begin{align*}
d_{\infty}(F(u), t'F(v))
&\geq d_{H}([F(u)]_{a}, [t'F(v)]_{t'^{-1}(a)})
=d_{H}(F([u]_{a}), F([v]_{t'^{-1}(a)}))\\
&=d_{H}\left(f_1([a, 1])\cup f_{2}([a, 1]), f_1\left(\left[a+\frac{1}{2}, 1\right]\right)\cup f_{2}\left(\left[a+\frac{1}{2}, 1\right]\right)\right)\\
&=d_{H}\left(\left[\frac{3}{16}, \frac{3}{4}\right], \left[\frac{7}{16}, \frac{1}{2}\right]\cup \left[\frac{21}{32}, \frac{3}{4}\right]\right)
=\frac{1}{4},
\end{align*}
implying that
$$
d_{0}(F(u), F(v))\geq \frac{1}{4}=d_{0}(u, v).
$$
Therefore, $F: (\mathfrak{F}(I), d_{0})\to (\mathfrak{F}(I), d_0)$ is not a contraction. This gives a negative
answer to Question~\ref{Q-2} as well.
\end{enumerate}
}
\end{remark}

\end{document}